\documentclass[12pt,reqno]{amsart}
\usepackage{amsmath, amsfonts, amssymb, amsthm,mathrsfs}
\def\Z{\mathbb Z}

\def\1{{\bf 1}}

\theoremstyle{plain}
\newtheorem{theorem}{Theorem}

\theoremstyle{definition}
\newtheorem*{acknowledgment}{Acknowledgments}
\theoremstyle{remark}

\begin{document}

\title{The restricted sumsets in finite abelian groups}
\author{Shanshan Du}
\address{ Department of Mathematics, Jinling Institute of Technology,
Nanjing 211169, People's Republic of China}
\email{ssdu@jit.edu.cn}

\author{Hao Pan}
\address{School of Applied Mathematics, Nanjing University of Finance and Economics, Nanjing 210046, People's Republic of China}
\email{haopan79@zoho.com}

\thanks{
The first author is supported by the National Natural Science Foundation of China (Grant No. 12371004, No. 11901259).
The second author is supported by the National Natural Science Foundation of China (Grant No. 12071208).}
\begin{abstract}
Suppose  that $k\geq 2$ and $A$ is a non-empty subset of a finite abelian group $G$ with $|G|>1$. Then the cardinality of the restricted sumset
$$
k^\wedge A:=\{a_1+\cdots+a_k:\,a_1,\ldots,a_k\in A,\ a_i\neq a_j\text{ for }i\neq j\}
$$
is at least
$$
\min\{p(G), k|A|-k^2+1\},
$$
where $p(G)$ denotes the least prime divisor of $|G|$.
\end{abstract}
\maketitle

\section{Introduction}\setcounter{equation}{0}
\setcounter{theorem}{0}
\setcounter{lemma}{0}
\setcounter{corollary}{0}
\setcounter{conjecture}{0}

Suppose that $G$ is a finite abelian group and $\emptyset\neq A,B\subseteq G$. Consider the sumset
$$
A+B:=\{a+b:\,a\in A,\ b\in B\}.
$$
For a prime $p$, let $\Z_p$ denote the cyclic group of order $p$. The classical Cauchy-Davenport theorem asserts that for any $\emptyset\neq A,B\subseteq\Z_p$,
\begin{equation}\label{CD}
|A+B|\geq\min\{p,|A|+|B|-1\}.
\end{equation}
In fact, according to the well-known Kneser theorem, (\ref{CD}) can be extended to any finite abelian group $G$ with $|G|>1$ as follows:
\begin{equation}
|A+B|\geq\min\{p(G),|A|+|B|-1\}
\end{equation}
for any non-empty subsets $A,B$ of $G$, where $p(G)$ denotes the least prime divisor of $|G|$.
Using an easy induction, we can get
\begin{equation}\label{CDG}
|A_1+\cdots+A_k|\geq\min\{p(G),|A_1|+\cdots+|A_k|-k+1\}
\end{equation}
for any $\emptyset\neq A_1,\ldots,A_k\subseteq G$.

On the other hand, for any $\emptyset\neq A\subseteq G$, consider the restricted sumset
$$
k^\wedge A:=\{a_1+\cdots+a_k:\,a_1,\ldots,a_k\in A,\ a_i\neq a_j\text{ for }i\neq j\}.
$$
In particular, $1^\wedge A=A$ and $0^\wedge A=\{0\}$.

In \cite{EH64}, Erd\H os and Heilbronn conjectured that for any prime $p$ and $\emptyset\neq A\subseteq\Z_p$,
\begin{equation}
|2^\wedge A|\geq\min\{p,2|A|-3\}.
\end{equation}
This conjecture was solved by Dias da Silva and Hamidoune \cite{DH94} with the help of the exterior algebra. In general they proved that
\begin{equation}\label{SHkA}
|k^\wedge A|\geq\min\{p,k|A|-k^2+1\}
\end{equation}
for any $\emptyset\neq A\subseteq\Z_p$. See \cite{ANR96} for an alternative proof based on the polynomial method.

In \cite{Ka04}, K\'arolyi considered the generalization of the Erd\H os-Heilbronn conjecture for finite abelian groups. Combining the polynomial method and a combinatorial induction, K\'arolyi obtained that
\begin{equation}\label{AAG}
|2^\wedge A|\geq\min\{p(G), 2|A|-3\}
\end{equation}
for any non-empty subset $A$ of a finite abelian group $G$ with $|G|>1$. An alternative proof of (\ref{AAG}) was given in \cite{Ka03}. For further
developments on restricted sumsets, see \cite{DP23}, \cite{KNPV15} and the references therein.

The purpose of our paper is to generalize (\ref{SHkA}) in the same spirit.

\begin{theorem}\label{mainT} Suppose that $G$ is a finite abelian group with $|G|>1$. For any $\emptyset\neq A\subseteq G$,
\begin{equation}\label{main}
|k^\wedge A|\geq\min\{p(G), k|A|-k^2+1\}.
\end{equation}
\end{theorem}

\bigskip
\section{The proof of Theorem \ref{mainT}}
\setcounter{equation}{0}
\setcounter{theorem}{0}
\setcounter{lemma}{0}
\setcounter{corollary}{0}
\setcounter{conjecture}{0}

We prove Theorem \ref{mainT} by induction on $|G|$. According to (\ref{SHkA}), it holds whenever $|G|$ is a prime,
so we suppose that $|G|$ is composite, and the result has been proved for smaller values of $|G|$. The result is trivial for $k=0,1$
or if $|A|=k$. It is also covered for $k=2$ by (\ref{AAG}), so in the sequel we will assume that $|A|\geq k+1\geq 4$.

Let $H$ be a subgroup of $G$ with $[G:H]=p(G)$. For each $a\in G$, consider the coset $\bar{a}=a+H.$ Let
$\bar{A}=\{\bar{a}:\,a\in A\}$, $m=|\bar{A}|$, and write
$$
A=\bigcup_{i=1}^m(a_i+A_i),
$$
where $a_1,\ldots,a_m\in A$ and $A_i\subseteq H$ with $a_{i}-a_{j}\not\in H$ for any distinct $i,j$. Without loss of generality, assume that
$
|A_1|\geq |A_2|\geq\cdots\geq |A_m|.
$
Note that $
|k^{\wedge}A|\geq |k^{\wedge}\bar{A}|$, so below assume that $k^{\wedge}\bar{A}\subsetneq G/H.$
Since (\ref{main}) follows from the induction hypothesis on $k^{\wedge}\bar{A}$ immediately when $|A|=m$, we will also assume that $|A_1|\geq 2.$

First we consider the relatively easy case when $m\geq k+1$.
For each $1\leq j\leq k+1$, let
$$
b_j=\sum_{\substack{1\leq i\leq k+1\\ i\neq j}}a_i
$$
and let
$$
U_j=A_1+\cdots+A_{j-1}+A_{j+1}+\cdots+A_{k+1}.
$$
Note that $\bar{b}_1,\ldots,\bar{b}_{k+1}$ are distinct elements of $G/H$. Hence
$$b_1+U_1,\ldots,b_{k+1}+U_{k+1}$$
are disjoint subsets of ${k^\wedge}A$.
Additionally, for each $1\leq j\leq k+1$, by (\ref{CDG}) we have
\begin{equation}\label{Bj+Ui}
|b_j+U_j|=\bigg|\sum_{\substack{1\leq i\leq k+1\\ i\neq j}}A_i\bigg|\geq\min\bigg\{p(H),
\sum_{\substack{1\leq i\leq k+1\\ i\neq j}}|A_i|-k+1\bigg\}.
\end{equation}

Furthermore, when $m\geq k+2$, for any $1\leq s\leq m-k-1$, let
$$
V_s=(k^\wedge\{\bar{a}_1,\ldots,\bar{a}_{k+s+1}\})\setminus(k^\wedge\{\bar{a}_1,\ldots,\bar{a}_{k+s}\}).
$$
Obviously,
$$
V_1=(k^\wedge\{\bar{a}_1,\ldots,\bar{a}_{k+2}\})\setminus\{\bar{b}_1,\ldots,\bar{b}_{k+1}\}.
$$
By (\ref{SHkA}) and the assumption that $k^{\wedge}\bar{A}\subsetneq G/H$, we have
\begin{align*}
\sum_{j=1}^s|V_j|=&|(k^\wedge\{\bar{a}_1,\ldots,\bar{a}_{k+s+1}\})\setminus\{\bar{b}_1,\ldots,\bar{b}_{k+1}\}|\\
\geq&k(k+s+1)-k^2+1-(k+1)=ks
\end{align*}
for each $1\leq s\leq m-k-1$.
On the other hand, if $\bar{c}\in V_s$, then $\bar{c}=\bar{a}_{i_1}+\cdots+\bar{a}_{i_k}$ where $1\leq i_1<\cdots<i_k\leq k+s+1$. It follows that
$$
|(k^\wedge A)\cap(\bar{c}+H)|\geq |A_{i_1}+\cdots+A_{i_k}|\geq|A_{i_1}|\geq|A_{s+2}|
$$
since $i_1\leq s+2$.
Note that if $x_1\geq\cdots\geq x_n\geq 0$ and the real numbers $\alpha_1, \ldots,\alpha_n$ satisfy $\sum_{i=1}^{r}\alpha_i\geq r$ for every $1\leq r\leq n$,
then $$\sum_{i=1}^{n}\alpha_i x_i\geq \sum_{i=1}^{n}x_i.$$
Hence we have
\begin{equation}\label{1mgeqk+2}
\sum_{j=1}^{m-k-1}|(k^\wedge A)\cap(V_j+H)|\geq\sum_{j=1}^{m-k-1}|V_j|\cdot|A_{j+2}|\geq k\sum_{j=1}^{m-k-1}|A_{j+2}|.
\end{equation}

In summary, using (\ref{Bj+Ui}) and (\ref{1mgeqk+2}),  we obtain that
\begin{align*}
&|k^\wedge A|\geq\sum_{j=1}^{k+1}|b_j+U_j|+\mathbf{1}_{m\geq k+2}\cdot\sum_{j=1}^{m-k-1}|(k^\wedge A)\cap(V_j+H)|\\
\geq&\min\bigg\{p(H),\sum_{j=1}^{k+1}\bigg(\sum_{\substack{1\leq i\leq k+1\\ i\neq j}}|A_i|-k+1\bigg)+\mathbf{1}_{m\geq k+2}\cdot\bigg(k\sum_{j=1}^{m-k-1} |A_{j+2}|\bigg)\bigg\}\\
\geq&\min\bigg\{p(H),k\sum_{j=1}^{k+1}|A_j|-k^2+1+\mathbf{1}_{m\geq k+2}\cdot\bigg(k\sum_{j=k+2}^{m}|A_j|\bigg)\bigg\}\\
\geq&\min\{p(G),k|A|-k^2+1\},
\end{align*}
where the indicative function $\mathbf{1}_{m\geq k+2}=\begin{cases} 1,&\text{if }m\geq k+2,\\
0,&\text{if }m=k+1.
\end{cases}$

\medskip
It remains to study the case $m\leq k$, which is more technical. Without loss of generality, assume that
$$
|A_1|=\cdots=|A_{t_1}|>|A_{t_1+1}|=\cdots=|A_{t_2}|>|A_{t_2+1}|=\cdots>|A_{t_{l-1}+1}|=\cdots=|A_{t_l}|,
$$
where $1\leq t_1<t_2<\cdots<t_l=m$. In addition, set $t_0=0$. We start with constructing some integers $n_1,\ldots,n_m\geq 1$ and $r$, $t$ with $1\leq r\leq t\leq m$ as follows.

Consider the sequence
\begin{align*}
&a_1,a_2,\ldots,a_m,\underbrace{a_{t_1},a_{t_1-1},\ldots,a_1,a_{t_1},\ldots,a_1,\ldots,a_{t_1},\ldots,a_1}_{(|A_{t_1}|-|A_{t_2}|)\text{ terms }a_{t_1},\ldots,a_1},\\
&\underbrace{a_{t_2},a_{t_2-1},\ldots,a_{1},a_{t_2},\ldots,a_{1},\ldots,a_{t_2},\ldots,a_{1}}_{(|A_{t_2}|-|A_{t_3}|)\text{ terms }a_{t_2},\ldots,a_{1}},\\
&\cdots\cdots,\\
&\underbrace{a_{t_l},a_{t_l-1},\ldots,a_{1},a_{t_l},\ldots,a_{1},\ldots,a_{t_l},\ldots,a_{1}}_{(|A_{t_l}|-1)\text{ terms }a_{t_l},\ldots,a_{1}}.
\end{align*}
Note that the length of the above sequence is
$$
m+\sum_{i=1}^{l-1} t_i\cdot (|A_{t_i}|-|A_{t_{i+1}}|)+t_l\cdot(|A_{t_l}|-1)=
t_1|A_{t_1}|+\sum_{i=2}^l(t_i-t_{i-1})\cdot|A_{t_i}|=|A|\geq k+1.
$$

Let $a_1^*,\ldots,a_{k+1}^*$ be the first $k+1$ terms of the above sequence,
then $a_1^*,\ldots,a_{k+1}^*$ must be of the form
\begin{align*}
&a_1,a_2,\ldots,a_m,a_{t_1},\ldots,a_1,\ldots,a_{t_1},\ldots,a_1,\\
&\cdots\cdots,\\
&a_{t_{s-1}},\ldots,a_{1},\ldots,a_{t_{s-1}},\ldots,a_{1},\\
&\underbrace{a_{t_s},\ldots,a_{1},\ldots,a_{t_s},\ldots,a_{1}}_{h\text{ terms }a_{t_s},\ldots,a_{1}},a_{t_s},\ldots,a_{r},
\end{align*}
where $1\leq r\leq t_s$ for some $1\leq s\leq l$, $0\leq h<|A_{t_s}|-|A_{t_{s+1}}|$ if $1\leq s\leq l-1$, and $0\leq h<|A_{t_s}|-1$ if $s=l$.
Now set $t:=t_s$, then we have $1\leq r\leq t$ and
\begin{equation}\label{HleqAt-2}
h\leq |A_t|-2.
\end{equation}

For each $1\leq i\leq m$, let
$$
n_i=|\{1\leq j\leq k+1:\,a_j^*=a_i\}|
.$$
Clearly we have
\begin{equation}\label{sumni}
\sum_{i=1}^mn_i=k+1,
\end{equation}
and when $t<i\leq m$, we have
\begin{equation}\label{ni1}
n_i=1.
\end{equation}
For any $i\in(t_{u-1},t_u]$ with $1\leq u\leq s$, if $1\leq i<r$, then
\begin{equation*}
n_i=1+h+\sum_{j=u}^{s-1}(|A_{t_j}|-|A_{t_{j+1}}|)=1+h+|A_{t_u}|-|A_{t_s}|=1+h+|A_i|-|A_{t}|,
\end{equation*}
and if $r\leq i\leq t$, then
\begin{equation*}
n_i=2+h+|A_i|-|A_{t}|.
\end{equation*}
That is,
\begin{equation}\label{Aini}
|A_i|-n_i=\begin{cases} |A_{t}|-h-1,&\text{if }1\leq i<r,\\
|A_{t}|-h-2,&\text{if }r\leq i\leq t.
\end{cases}
\end{equation}
Then by (\ref{HleqAt-2}), (\ref{ni1}) and (\ref{Aini}), we have
\begin{equation}\label{nileqAi}
n_i\leq |A_i|
\text{ for each } r\leq i\leq m,
\end{equation}
\begin{equation}\label{nileqAi-1}
n_i\leq |A_i|-1
\text{ for each } 1\leq i< r,
\end{equation}
and \begin{equation}\label{AiniAjnj1}
|A_i|-n_i=|A_j|-n_j+1
\text{ for each } 1\leq i< r \text{ and } r\leq j\leq t.
\end{equation}
Furthermore, for $1\leq i\leq t$ and $t<j\leq m$, we have
\begin{equation}\label{AiniAjnj}
|A_i|-n_i\geq|A_{t}|-h-2=|A_{t_s}|-h-2\geq |A_{t_{s+1}}|-1\geq|A_j|-n_j,
\end{equation}
since now $s<l$.

For non-negative integers $l_1,\ldots,l_m$, let
$$
S_{l_1,\ldots,l_m}={l_1}^\wedge A_1+{l_2}^\wedge A_2+\cdots+{l_m}^\wedge A_m\subseteq H.
$$
For each $1\leq j\leq m$, let
$$b_j=(n_j-1)a_j+\sum_{\substack{1\leq i\leq m\\ i\neq j}}n_ia_i,
$$
and let
$$
U_j=S_{n_1,\ldots,n_{j-1},n_j-1,n_{j+1},\ldots,n_m}.
$$
Since $n_i\leq |A_i|$, we know that $U_j\neq\emptyset$.
Clearly $b_1+U_1,\ldots,b_m+U_m$ are disjoint subsets of $k^\wedge A$. According to the induction hypothesis, we have
$$
|n^\wedge A_i|\geq\min\{p(H),n|A_i|-n^2+1\}
$$
for any $n\geq 0$ and $1\leq i\leq m$.
Hence by (\ref{CDG}),
\begin{align}\label{Uj}
&|U_j|=\bigg|(n_j-1)^\wedge A_j+\sum_{\substack{1\leq i\leq m\\ i\neq j}}{n_i}^\wedge A_i\bigg|\notag\\
\geq&\min\bigg\{p(H),(n_j-1)|A_j|-((n_j-1)^2-1)+\sum_{\substack{1\leq i\leq m\\ i\neq j}}n_i|A_i|-\sum_{\substack{1\leq i\leq m\\ i\neq j}}(n_i^2-1)-m+1\bigg\}\notag\\
=&\min\bigg\{p(H),(n_j-1)|A_j|+\sum_{\substack{1\leq i\leq m\\ i\neq j}}n_i|A_i|-\sum_{i=1}^mn_i^2+2n_j\bigg\}.
\end{align}

However, we need to find more elements lying in $k^\wedge A$. For any $I\subseteq\{1,2,\ldots,t\}$ with $|I|=t-r+1$,
we have
\begin{equation}\label{bjInj}
b_j-\sum_{i=r}^{t}a_i+\sum_{i\in I}a_i=(n_{j,I}-1)a_j+\sum_{\substack{1\leq i\leq m\\ i\neq j}}n_{i,I} a_i
\end{equation}
for $1\leq j\leq m$, where
$$
n_{i,I}=\begin{cases}
n_i-1,&\text{if }i\in[r,t]\setminus I,\\
n_i+1,&\text{if }i\in I\setminus [r,t],\\
n_i,&\text{otherwise}
\end{cases}
$$
for each $1\leq i\leq m$.
If $r\leq  i\leq t$, then $n_i\geq 2$. So we always have $n_{i,I}\geq 1$. Furthermore, by (\ref{nileqAi}) and (\ref{nileqAi-1}), we obtain that $n_{i,I}\leq |A_i|$ for each $1\leq i\leq m$. Then \begin{equation}\label{SnI}
S_{n_{1,I},\ldots,n_{j-1,I},n_{j,I}-1,n_{j+1,I},\ldots,n_{m,I}}\neq\emptyset
\end{equation}
for any $1\leq j\leq m$.

Let
$$
d_j=b_j-\sum_{i=r}^{t}a_i
$$
for $1\leq j\leq m$, and let
$$
R=\{\bar{d}_1,\ldots,\bar{d}_m\}+(t-r+1)^\wedge\{\bar{a}_1,\bar{a}_2,\ldots,\bar{a}_{t}\}\subseteq G/H.
$$
Then
\begin{align}
|R|&\geq\min\big\{p(G), m+((t-r+1)\cdot t-(t-r+1)^2+1)-1\big\}\notag\\
&=\min\big\{p(G),m+(t-r+1)\cdot(r-1)\big\}.\notag
\end{align}
First we suppose $m+(t-r+1)(r-1)<p(G)$. Obviously $(t-r+1)(r-1)\geq 0$. If $(t-r+1)(r-1)\geq 1$, then there exist $$1\leq j_1,\ldots,j_{(t-r+1)(r-1)}\leq m$$ and $$I_1,\ldots,I_{(t-r+1)(r-1)}\subseteq\{1,\ldots,t\}$$ with $|I_j|=t-r+1$ for
$1\leq j\leq (t-r+1)(r-1)$
such that the elements
$$
\bar{c}_1,\bar{c}_2,\ldots,\bar{c}_{(t-r+1)(r-1)}\in R\setminus\{\bar{b}_1,\ldots,\bar{b}_m\}
$$
are distinct, where
$$
c_u=d_{j_u}+\sum_{i\in I_{u}}a_i
$$
for each $1\leq u\leq (t-r+1)(r-1)$.

Let
$$
V_u=S_{n_{1,I_u},\ldots,n_{j_{u}-1,I_u},n_{j_u,I_u}-1,n_{j_u+1,I_u},\ldots,n_{m,I_u}}
$$
for $1\leq u\leq (t-r+1)(r-1)$. Thus in view of (\ref{bjInj}) and (\ref{SnI}), we know that
$$
b_1+U_1,\ldots,b_{m}+U_m,c_1+V_1,\ldots,c_{(t-r+1)(r-1)}+V_{(t-r+1)(r-1)}
$$
are disjoint non-empty subsets of $k^\wedge A$.

It follows from (\ref{Uj}) and (\ref{SnI}) that
\begin{align}\label{kAUj}
&|k^\wedge A|\geq\sum_{j=1}^m|U_j|+(t-r+1)(r-1)\notag\\
\geq&\min\bigg\{p(H),\sum_{j=1}^m\bigg((n_j-1)|A_j|+\sum_{\substack{1\leq i\leq m\\ i\neq j}}n_i|A_i|-\sum_{i=1}^mn_i^2+2n_j\bigg)+(t-r+1)(r-1)\bigg\}\notag\\
=&\min\bigg\{p(H),m\sum_{i=1}^mn_i|A_i|-m\sum_{i=1}^mn_i^2+2\sum_{i=1}^mn_i-|A|+(t-r+1)(r-1)\bigg\}.
\end{align}
Note that
\begin{align*}
m\sum_{i=1}^mn_i|A_i|=&\bigg(\sum_{i=1}^mn_i\bigg)\cdot\bigg(\sum_{i=1}^m|A_i|\bigg)+\sum_{\substack{1\leq i,j\leq m\\ i\neq j}}n_i(|A_i|-|A_j|).
\end{align*}
So by (\ref{sumni}), we have
\begin{align}\label{mniAi}
&m\sum_{i=1}^mn_i|A_i|-m\sum_{i=1}^mn_i^2+2\sum_{i=1}^mn_i-|A|\notag\\
=&k|A|+\sum_{\substack{1\leq i,j\leq m\\ i\neq j}}n_i(|A_i|-|A_j|)-k^2+\bigg(\sum_{i=1}^mn_i-1\bigg)^2-m\sum_{i=1}^mn_i^2+2\sum_{i=1}^mn_i\notag\\
=&k|A|-k^2+1+\sum_{\substack{1\leq i,j\leq m\\ j>i}}(n_i-n_j)\cdot(|A_i|-|A_j|)+\bigg(\sum_{i=1}^mn_i\bigg)^2-m\sum_{i=1}^mn_i^2.
\end{align}

Furthermore by (\ref{ni1}), (\ref{Aini}), (\ref{AiniAjnj1}) and (\ref{AiniAjnj}), we have
$$
|A_i|-|A_j|\ \ \begin{cases}\geq n_i-n_j\geq 0,&\text{if }1\leq i\leq t<j\leq m\text{, or }t<i<j\leq m,\\
=n_i-n_j\geq 0,&\text{if }1\leq i<j<r\text{, or }r\leq i<j\leq t,\\
=n_i-n_j+1\geq 0,&\text{if }1\leq i<r\leq j\leq t.
\end{cases}
$$
Hence for any $1\leq i<j\leq m$, if $1\leq i<r\leq j\leq t$, we have
$$
(n_i-n_j)\cdot(|A_i|-|A_j|)= (n_i-n_j)^2+(n_i-n_j)
$$
and $n_i-n_j\geq -1$, otherwise we have
$$
(n_i-n_j)\cdot(|A_i|-|A_j|)\geq (n_i-n_j)^2.
$$
 It follows that
\begin{align}\label{ninjAiAj}
&\sum_{\substack{1\leq i,j\leq m\\ j>i}}(n_i-n_j)(|A_i|-|A_j|)+\bigg(\sum_{i=1}^mn_i\bigg)^2-m\sum_{i=1}^mn_i^2\notag\\
\geq&\sum_{\substack{1\leq i,j\leq m\\ j>i}}(n_i-n_j)^2+\sum_{\substack{1\leq i<r\\ r\leq j\leq t}}(n_i-n_j)+\bigg(\sum_{i=1}^mn_i\bigg)^2-m\sum_{i=1}^mn_i^2\notag\\
\geq&(m-1)\sum_{i=1}^mn_i^2-2\sum_{\substack{1\leq i,j\leq m\\ j>i}}n_in_j-\sum_{\substack{1\leq i<r\\ r\leq j\leq t}}1+\sum_{i=1}^mn_i^2+2\sum_{\substack{1\leq i,j\leq m\\ j>i}}n_in_j-m\sum_{i=1}^mn_i^2\notag\\
=&-(r-1)\cdot(t-r+1).
\end{align}

Combining (\ref{kAUj}), (\ref{mniAi}) and (\ref{ninjAiAj}), we obtain that
$$
|k^\wedge A|\geq
\min\{p(H),k|A|-k^2+1\}\geq
\min\{p(G),k|A|-k^2+1\}.
$$

Finally, if $m+(t-r+1)(r-1)\geq p(G)$, then by a similar discussion,
one can deduce that
$$|k^\wedge A|\geq |R|\geq p(G).$$

\qed

\begin{acknowledgment}
We are grateful to the anonymous referee for his/her very useful comments and helpful suggestions on our paper.
\end{acknowledgment}

\end{document}